\documentclass[12pt,a4paper,dvips]{article}
\usepackage{inputenc, amssymb, amsmath, amsthm, srcltx}
\usepackage[russian,english]{babel}
\usepackage{color}
\usepackage{graphicx}
\makeatletter
\thm@style{plain}

\def\th@plain{%
  \thm@headfont{\bfseries}%
  \itshape 
  \thm@notefont{\rm}%
}
\def\thm@indent{\hspace*{\parindent}}
\makeatother

\def\({\left(}
\def\){\right)}

\newcommand{\be}{\begin{equation}}
\newcommand{\ee}{\end{equation}}

\renewcommand{\geq}{\geqslant}
\let\epsilon\varepsilon
\let\phi\varphi
\let\le\leqslant
\let\ge\geqslant

\def\arg{\mathop{\ensuremath{\text{\textup{arg}}}}}
\def\sgntwods{\mathop{\ensuremath{\text{\textup{sgn}}_{2\Delta\sigma}}}}
\def\sgn{\mathop{\ensuremath{\text{\textup{sgn}}}}}
\def\Re{\mathop{\ensuremath{\text{\textup{Re}}}}}
\def\Im{\mathop{\ensuremath{\text{\textup{Im}}}}}

\newtheorem{theorem}{Theorem}
\newtheorem{lemma}{Lemma}
\newcommand{\dokvo}{{\it Proof.} }
\newcommand\inte{\int\limits}

\newcommand{\beq}{\begin{equation}}
\newcommand{\eeq}{\end{equation}}

{Definition}
{Algorithm}

\newcommand{\ds}{\Delta\sigma}
\newcommand{\K}{\mathcal{K}}
\newcommand{\RR}{\mathcal{R}}

\begin{document}

\centerline{\bf\uppercase{On a choice of the mollified function}}
\centerline{\bf\uppercase{in the Levinson--Conrey method}\footnote[1]{%
}}

\bigskip

\medskip

\centerline{\sc Sergei~Preobrazhenski\u i and Tatyana~Preobrazhenskaya}

\bigskip

\bigskip

\hbox to \textwidth{\hfil\parbox{0.9\textwidth}{%
\small {\sc Abstract.}
Motivated by a functional property of the Riemann zeta function,
we consider a new form of the mollified function
in the Levinson--Conrey method.
As an application, we give the following slight improvement of Feng's result:
assuming Feng's condition on the lengths of the mollifier at least $41{.}2948$\% of the zeros
of the Riemann zeta function are on the critical line.

The construction may lead to further improvements as one increases the number of terms
in Feng's mollifier.}\hfil}

\bigskip

\bigskip
\textbf{Contents}\hfil

\bigskip

1. Introduction

\smallskip

2. Main lemma

\smallskip

3. Proof of Theorem 1

\smallskip

4. Further remarks

\smallskip

5. Appendix

\smallskip

References

\bigskip

\bigskip

\bigskip


{\bf 1. Introduction.} The Riemann zeta-function $\zeta(s)$ is defined for $\Re s>1$ by
\[
\zeta(s)=\sum_{n=1}^{\infty}n^{-s},
\]
and for other $s$ by the analytic continuation.
It is a meromorphic function in the whole complex plane
with the only singularity $s=1$, which is a simple pole
with residue $1$.

The Euler product links the zeta-function
and prime numbers: for $\Re s>1$
\[
\zeta(s)=\prod_{p\text{\textup{ prime}}}\left(1-p^{-s}\right)^{-1}.
\]

The functional equation for $\zeta(s)$ may be written in the form
\[
\xi(s)=\xi(1-s),
\]
where $\xi(s)$ is an entire function defined by
\[
\xi(s)=H(s)\zeta(s)
\]
with
\[
H(s)=\frac12s(1-s)\pi^{-s/2}\Gamma\left(\frac s2\right).
\]
This implies that $\zeta(s)$ has zeros at
$s=-2$, $-4$, ${\ldots}$
These zeros are called the ``trivial'' zeros.
It is known that $\zeta(s)$ has infinitely many
nontrivial zeros $s=\rho=\beta+i\gamma$, and all of them are in the ``critical strip''
$0<\Re s=\sigma<1$, $-\infty<\Im s=t<\infty$.
The pair of nontrivial zeros with the smallest value of $|\gamma|$
is $\frac12\pm i(14{.}134725\ldots)$.

If $N(T)$ denotes the number of zeros $\rho=\beta+i\gamma$
($\beta$ and $\gamma$ real), for which $0<\gamma\le T$,
then
\[
N(T)=\frac T{2\pi}\log\left(\frac T{2\pi}\right)
-\frac T{2\pi}+\frac78+S(T)+O\left(\frac1T\right),
\]
with
\[
S(T)=\frac1{\pi}\arg\zeta\left(\frac12+iT\right)
\]
and
\[
S(T)=O(\log T).
\]
This is the Riemann--von~Mangoldt formula for $N(T)$.

Let $N_0(T)$ be the number of zeros of $\zeta\left(\frac12+it\right)$ when $0<t\le T$,
each zero counted with multiplicity.
The Riemann hypothesis is the conjecture that $N_0(T)=N(T)$.
Let
\[
\kappa=\liminf\limits_{T\to\infty}\frac{N_0(T)}{N(T)}.
\]

Important results about $N_0(T)$ include:
\begin{itemize}
\item \cite{PreobHardy14}: Hardy proved that $N_0(T)\to\infty$ as $T\to\infty$.
\item \cite{PreobHL21}: Hardy and Littlewood obtained that $N_0(T)\ge AT$ for some $A>0$
and all sufficiently large $T$.
\item \cite{PreobSel42}: Selberg proved that $\kappa\ge A$ for an effectively computable positive constant $A$.
\item \cite{PreobLev74}: Levinson proved that $\kappa\ge0{.}34\ldots$
\item \cite{PreobCon89}: Conrey obtained $\kappa\ge0{.}4088\ldots$
\item \cite{PreobFeng12}: Feng obtained $\kappa\ge0{.}4128\ldots$ (assuming a condition on the lengths of the mollifier)
\end{itemize}

In this article we establish the following improvement of Feng's result:
\begin{theorem}\label{PreobLargeProportion}
If Feng's asymptotic formula for the mean square,
used to estimate the integral $I(R)$ in Theorem~\textup{\ref{PreobLevCon}}, is valid
for the mollifier $M(s)$ with $\theta=\frac47-\varepsilon$ and $\theta_1=\frac12-\varepsilon$ then we have
\[
\kappa\ge0{.}412948\ldots
\]
\end{theorem}

The motivation for our choice of the mollified function is Lemma~\ref{PreobConreyShLemma}
which allows to ``translate'' certain terms of the sum
in general Conrey's construction~\cite{PreobCon83} by
\[
\Delta\sigma=\frac{\alpha}{\log T}.
\]
It turns out that the translating is relevant when we increase the number of terms
in Feng's mollifier.

\textbf{2. Main lemma.}
\begin{lemma}\label{PreobConreyShLemma}
Let $f(s)$ be an analytic function, $s\in{\mathbb C}$, $\ds\in{\mathbb R}$,
$\K\ge1$ be an odd integer.
Then
\be\label{PreobShIdentity}
\begin{split}
&f(s+\ds)=f(s)+\sum_{\substack{k\text{\textup{ odd}}\\ k\le\K}}\left(g_k(\ds)f^{(k)}(s)+g_k(\ds)f^{(k)}(s+\ds)\right)\\
&+\frac{4(-1)^{(\K+1)/2}(\ds)^{\K+1}}{\pi^{\K+2}}\inte_s^{s+\ds}f^{(\K+2)}(w)\left(\sum_{n=1}^{\infty}\frac1{(2n-1)^{\K+2}}\sin\left(\frac{(2n-1)\pi(s+\ds-w)}{\ds}\right)\right)\,dw,
\end{split}
\ee
where
\[
g_k(\ds)=\frac{4(-1)^{(k-1)/2}(\ds)^k}{\pi^{k+1}}\sum_{n=1}^{\infty}\frac1{(2n-1)^{k+1}}.
\]
\end{lemma}

\dokvo We use induction on $\K$. First establish induction base $\K=1$. We have
\[
f(s+\ds)=f(s)+\inte_s^{s+\ds}f'(w)\,dw.
\]
Let $\sgntwods(x)$ be the $2\ds$-periodic real-valued function defined by
\[
\sgntwods(x)=\begin{cases}
1&\text{if $x\in(0,\ds)$},\\
0&\text{if $x=-\ds,0,\ds$},\\
-1&\text{if $x\in(-\ds,0)$}.
\end{cases}
\]
Using the Fourier expansion
\be\label{Preobsgnfe}
\sgntwods(x)=\frac4{\pi}\sum_{n=1}^{\infty}\frac1{2n-1}\sin\left(\frac{(2n-1)\pi x}{\ds}\right)
\ee
we obtain
\[
\begin{split}
f(s+\ds)&=f(s)+\frac4{\pi}\inte_s^{s+\ds}f'(w)\left(\sum_{n=1}^{\infty}\frac1{2n-1}\sin\left(\frac{(2n-1)\pi(s+\ds-w)}{\ds}\right)\right)\,dw,\\
f(s+\ds)&=f(s)+\frac4{\pi}\left(\inte_s^{s+\epsilon}\cdots+\inte_{s+\epsilon}^{s+\ds-\epsilon}\cdots+\inte_{s+\ds-\epsilon}^{s+\ds}\cdots\right).
\end{split}
\]
The series~\eqref{Preobsgnfe} converges uniformly in $x\in[\epsilon,\ds-\epsilon]$ so by integrating by parts
\[
\begin{split}
&\frac4{\pi}\inte_{s+\epsilon}^{s+\ds-\epsilon}f'(w)\,d\left(\sum_{n=1}^{\infty}\frac{\ds}{(2n-1)^2\pi}\cos\left(\frac{(2n-1)\pi(s+\ds-w)}{\ds}\right)\right)\\
\mathrel{=}&\frac{4\ds}{\pi^2}\sum_{n=1}^{\infty}\frac1{(2n-1)^2}\Bigl(f'(s+\ds-\epsilon)+f'(s+\epsilon)\Bigr)\\
-&\frac{4\ds}{\pi^2}\inte_{s+\epsilon}^{s+\ds-\epsilon}f''(w)\left(\sum_{n=1}^{\infty}\frac1{(2n-1)^2}\cos\left(\frac{(2n-1)\pi(s+\ds-w)}{\ds}\right)\right)\,dw+\delta_1(\epsilon)\\
\mathrel{=}&g_1(\ds)\Bigl(f'(s)+f'(s+\ds)\Bigr)\\
-&\frac{4\ds}{\pi^2}\inte_s^{s+\ds}f''(w)\,d\left(\sum_{n=1}^{\infty}-\frac{\ds}{(2n-1)^3\pi}\sin\left(\frac{(2n-1)\pi(s+\ds-w)}{\ds}\right)\right)+\delta_2(\epsilon)\\
\mathrel{=}&g_1(\ds)\Bigl(f'(s)+f'(s+\ds)\Bigr)\\
-&\frac{4\ds}{\pi^2}\left(-\inte_s^{s+\ds}-\frac{\ds}{\pi}f'''(w)\left(\sum_{n=1}^{\infty}\frac1{(2n-1)^3}\sin\left(\frac{(2n-1)\pi(s+\ds-w)}{\ds}\right)\right)dw\right)+\delta_3(\epsilon),
\end{split}
\]
and $\delta_1(\epsilon),\delta_2(\epsilon),\delta_3(\epsilon)\to0$ as $\epsilon\to0$.
This proves the induction base. The induction step is proven by integrating by parts in~\eqref{PreobShIdentity}
as above,
with the uniform convergence of the series in the integrand when $\K\ge1$.

\noindent{\sc Remark.} We have
\[
g_1(\ds)=\frac{4\ds}{\pi^2}\sum_{n=1}^{\infty}\frac1{(2n-1)^2}=\frac{4\ds}{\pi^2}\zeta(2)\left(1-\frac1{2^2}\right)=\frac{\ds}2,
\]
and in general for $k$ odd
\be\label{PreobBernoulliEq}
g_k(\ds)=\frac{4(-1)^{(k-1)/2}(\ds)^k}{\pi^{k+1}}\zeta(k+1)\left(1-\frac1{2^{k+1}}\right)
=-\left(\frac{\ds}2\right)^k\frac{2^{k+1}-4^{k+1}}{(k+1)!}B_{k+1},
\ee
where $B_{k+1}$ is the Bernoulli number.

The series
\[
\sum_{\substack{k\ge1\\ k\text{\textup{ odd}}}}\left(g_k(\ds)f^{(k)}(s)+g_k(\ds)f^{(k)}(s+\ds)\right)
\]
obtained by successive integrations by parts in~\eqref{PreobShIdentity}
may be divergent. However, we have the following
\begin{lemma}\label{PreobBernoulliLemma}
Suppose that $0<\epsilon<2\pi$, $|\alpha|\le2\pi-\epsilon$
and
\[
|\ds|=\frac{|\alpha|}{\log T}\le\frac{2\pi-\epsilon}{\log T}.
\]
Then the series
\[
\sum_{\substack{k\ge1\\ k\text{\textup{ odd}}}}\Bigl(-g_k(\ds)\Bigr)(\log T)^k\left(\frac12-x\right)^k
\]
converges on $x\in[0,1]$ and
\[
\sum_{\substack{k\ge1\\ k\text{\textup{ odd}}}}\Bigl(-g_k(\ds)\Bigr)(\log T)^k\left(\frac12-x\right)^k
=-\tanh\left(\frac{\alpha}2\left(\frac12-x\right)\right). 
\]
\end{lemma}

\dokvo From~\eqref{PreobBernoulliEq} we have
\[
\Bigl(-g_k(\ds)\Bigr)(\log T)^k\left(\frac12-x\right)^k=\frac{2\left(\alpha\left(\frac12-x\right)\right)^kB_{k+1}}{(k+1)!}-\frac{4(\alpha(1-2x))^kB_{k+1}}{(k+1)!}.
\]
By the definition of the Bernoulli numbers,
\[
\frac{z}{e^z-1}=\sum_{m=0}^{\infty}\frac{B_m}{m!}z^m,
\]
with $B_0=1$, $B_1=-\frac12$ and $B_3=B_5=B_7=\dots=0$,
the radius of convergence of the series being $2\pi$.
Then
\[
\begin{split}
&2\sum_{\substack{k\ge1\\ k\text{\textup{ odd}}}}\frac{\left(\alpha\left(\frac12-x\right)\right)^kB_{k+1}}{(k+1)!}
-4\sum_{\substack{k\ge1\\ k\text{\textup{ odd}}}}\frac{(\alpha(1-2x))^kB_{k+1}}{(k+1)!}\\
=&\frac2{\alpha(1/2-x)}\left(\frac{\alpha(1/2-x)}{e^{\alpha(1/2-x)}-1}-1+\frac{\alpha(1/2-x)}2\right)
-\frac4{\alpha(1-2x)}\left(\frac{\alpha(1-2x)}{e^{\alpha(1-2x)}-1}-1+\frac{\alpha(1-2x)}2\right),
\end{split}
\]
and the lemma follows.

\begin{lemma}\label{PreobHsplusdsLemma}
Suppose that $\alpha$ is real
and
\[
\ds=\frac{\alpha}{\log T}.
\]
Then in the rectangle
\[
s=\sigma+it,\quad\frac13\le\sigma\le A,\quad T\le t\le2T
\]
with $A\ge3$ and $T\ge2A$
we have
\[
H(s+\ds)=\left(e^{\alpha/2}+O\left(\frac1{\log T}\right)\right)H(s).
\]
\end{lemma}

\textbf{3. Proof of Theorem~\ref{PreobLargeProportion}.}
Here we give a sketch of the argument and motivate our choice of the mollified function
(see Iwaniec' lecture notes~\cite{PreobIw14}).

Suppose that $0<\epsilon<2\pi$, $|\alpha|=2\pi-\epsilon$
and
\[
|\ds|=\frac{|\alpha|}{\log T}=\frac{2\pi-\epsilon}{\log T}.
\]
Define $G(s)=G_{\epsilon_1,\K,\ds}(s)$ by
\[
2e^{-\alpha/2}H(s+\ds)G(s)=\xi(s)+\epsilon_1\left(\tilde{g}_0\xi(s)+\sum_{\substack{k\text{\textup{ odd}}\\ k\le\K_0}}\tilde{g}_k\xi^{(k)}(s)\right)+\sum_{\substack{k\text{\textup{ odd}}\\ k\le\K}}g_k(\ds)\xi^{(k)}(s),
\]
where $\epsilon_1$ is a real number to be chosen later.
For $k$ odd $\xi^{(k)}(s)$ is purely imaginary on the line $\Re s=\frac12$, so
\[
\Re2e^{-\alpha/2}H(s+\ds)G(s)=(1+\epsilon_1\tilde{g}_0)\xi(s)\qquad\text{if }\Re s=\frac12.
\]
That is, the critical zeros of $\zeta(s)$ are precisely the points on the line $\Re s=\frac12$
for which we have either $G(s)=0$
or
\[
G(s)\ne0\qquad\text{and}\qquad\!\arg H(s+\ds)G(s)\equiv\frac{\pi}2\pmod\pi.
\]
Next for
\[
|\alpha|=2\pi-\epsilon,\qquad|\ds|=\frac{|\alpha|}{\log T}=\frac{2\pi-\epsilon}{\log T}
\]
we obtain (see~\cite{PreobIw14})
\[
G(s)=\sum_{l\le T}Q_{1,\epsilon_1,\K,\alpha}\left(\frac{\log l}{\log T}+\delta_1(s)\right)l^{-s}+O\left(T^{-\frac14}\right),
\]
where
\[
\delta_1(s)\ll\frac1{\log T}
\]
and $Q_{1,\epsilon_1,\K,\alpha}(x)$ is the polynomial
\[
Q_{1,\epsilon_1,\K,\alpha}(x)=\frac12+\epsilon_1Q_0(x)+\frac2{\pi}\sum_{\substack{k\text{\textup{ odd}}\\ k\le\K}}
(\sgn\alpha)^k\left(1-\frac{\epsilon}{2\pi}\right)^k(-1)^{(k-1)/2}\zeta(k+1)\left(1-\frac1{2^{k+1}}\right)(1-2x)^k.
\]
Now we choose $\epsilon_1$ so that the polynomial $Q_{1,\epsilon_1,\K,\alpha}(x)$ satisfies
\[
Q_{1,\epsilon_1,\K,\alpha}(0)=1.
\]
Let $N_{01}(T,2T)$ denote the number of zeros of $\zeta(s)$ when
\[
s=\rho=\frac12+i\gamma,\qquad T\le\gamma\le2T,
\]
counted without multiplicity.
By estimating the argument variations we get
\[
N_{01}(T,2T)\ge N(T,2T)-2N_G(\RR)+O(T),
\]
where $N(T,2T)$ is the number of all zeros $\rho=\beta+i\gamma$ of $\zeta(s)$ with
$T\le\gamma\le2T$ counted with multiplicity,
and $N_G(\RR)$ denotes the number of zeros counted with multiplicity of $G(s)$ inside
the closed rectangle $\RR$ that has the segment
\[
\Re s=\frac12,\qquad T\le\Im s\le2T
\]
as its left side,
with small circular dents to the right centered at the common critical zeros of $\zeta(s)$
and $G(s)$,
and that has the segment
\[
\Re s=A,\qquad T\le\Im s\le2T
\]
with $A$ a large enough constant as its right side.

By Lemma~\ref{PreobConreyShLemma} we can write
\[
\begin{split}
2e^{-\alpha/2}H(s+\ds)G(s)&=\xi(s+\ds)
+\sum_{\substack{k\text{\textup{ odd}}\\ k\le\K}}\Bigl(-g_k(\ds)\Bigr)\xi^{(k)}(s+\ds)\\
&+\epsilon_1\left(\tilde{g}_0\xi(s)+\sum_{\substack{k\text{\textup{ odd}}\\ k\le\K_0}}\tilde{g}_k\xi^{(k)}(s)\right)
+\RR_1(\K,s,\ds,\epsilon_1),
\end{split}
\]
where
\[
\begin{split}
&\RR_1(\K,s,\ds,\epsilon_1)\\
&=\frac{4(-1)^{(\K-1)/2}(\ds)^{\K+1}}{\pi^{\K+2}}\inte_s^{s+\ds}\xi^{(\K+2)}(w)
\left(\sum_{n=1}^{\infty}\frac1{(2n-1)^{\K+2}}\sin\left(\frac{(2n-1)\pi(s+\ds-w)}{\ds}\right)\right)\,dw,
\end{split}
\]
and
\[
-g_k(\ds)=\frac{4(-1)^{(k+1)/2}(\ds)^k}{\pi^{k+1}}\zeta(k+1)\left(1-\frac1{2^{k+1}}\right).
\]

From this we obtain another representation for $e^{-\alpha/2}G(s)=e^{-\alpha/2}G_{\epsilon_1,\K,\ds}(s)$:
\[
\begin{split}
e^{-\alpha/2}G(s)&=\sum_{l\le T}Q_{\K,\alpha}\left(\frac{\log l}{\log T}+\delta(s)\right)l^{-(s+\ds)}
+\epsilon_1e^{-\alpha/2}\sum_{l\le T}Q_0\left(\frac{\log l}{\log T}+\delta_0(s)\right)l^{-s}\\
&+\sum_{l\le T}Q_{2,\K,\ds}\left(\frac{\log l}{\log T}+\delta_{2,\K,\ds}(s)\right)l^{-(s+\ds)}
+\sum_{l\le T}Q_{3,\K,\ds}\left(\frac{\log l}{\log T}+\delta_{3,\K,\ds}(s)\right)l^{-s},
\end{split}
\]
where
\[
\delta(s),\delta_0(s),\delta_{2,\K,\ds}(s),\delta_{3,\K,\ds}(s)\ll\frac1{\log T}
\]
and $Q_{\K,\alpha}(x)$ is the polynomial
\[
Q_{\K,\alpha}(x)=\frac12+\frac2{\pi}\sum_{\substack{k\text{\textup{ odd}}\\ k\le\K}}
(\sgn\alpha)^k\left(1-\frac{\epsilon}{2\pi}\right)^k(-1)^{(k+1)/2}\zeta(k+1)\left(1-\frac1{2^{k+1}}\right)(1-2x)^k.
\]
The polynomials $Q_{2,\K,\ds}(x)$ and $Q_{3,\K,\ds}(x)$ are obtained
by making the integration in
\[
\RR_1(\K,s,\ds,\epsilon_1).
\]
Note that by Lemma~\ref{PreobBernoulliLemma} we have
\be\label{Q1x}
Q_1(x)=\lim\limits_{\substack{\K\to\infty\\ |\alpha|<2\pi}}Q_{1,\epsilon_1,\K,\alpha}(x)=\frac12-\frac12\tanh\left(\frac{\alpha}2\left(x-\frac12\right)\right)+\epsilon_1Q_0(x),
\ee
so that $Q_1(0)=1$ implies $\epsilon_1Q_0(0)=\frac12-\frac12\tanh\left(\frac{\alpha}4\right)$ and
\[
Q_1(x)=\frac12-\frac12\tanh\left(\frac{\alpha}2\left(x-\frac12\right)\right)+Q_0^{-1}(0)\left(\frac12-\frac12\tanh\left(\frac{\alpha}4\right)\right)Q_0(x).
\]
Also by Lemma~\ref{PreobBernoulliLemma}
\[
Q(x)=\lim\limits_{\substack{\K\to\infty\\ |\alpha|<2\pi}}Q_{\K,\alpha}(x)=\frac12+\frac12\tanh\left(\frac{\alpha}2\left(x-\frac12\right)\right),
\]
hence
\[
Q(0)=\frac12-\frac12\tanh\left(\frac{\alpha}4\right).
\]
Note that the term
\[
\frac12-\frac12\tanh\left(\frac{\alpha}2\left(x-\frac12\right)\right)
\]
in $Q_1(x)$ transforms to $Q(x)$, and the corresponding sum in the expression for $G(s)$
is now translated by $\ds$.
This motivates our choice of the function $Q_1(x)$ in the form
\[
Q_1(x)=\frac12-\frac12\tanh\left(\frac{\alpha}2\left(x-\frac12\right)\right)+\left(\frac12-\frac12\tanh\left(\frac{\alpha}4\right)\right)(2\tilde{Q}(x)-1)
\]
with some appropriate function $\tilde{Q}(x)$, e.g. a polynomial satisfying $\tilde{Q}(0)=1$ and $\tilde{Q}(x)+\tilde{Q}(1-x)\equiv g$ for some real $g$.

In fact, the polynomial $\tilde{Q}(x)$ will be taken the near-optimal polynomial of degree $d_n$
obtained by Feng's method in step $n$ using polynomials $P_{1,n}$, $\ldots$, $P_{I_n,n}$ in the mollifier.
Next, the $d_{n+1}$-truncation of the series for $Q_1(x)$ will be used in step $n+1$
to optimize the value of $\kappa$ over $\alpha=\alpha_{n+1}$ and the new polynomials $P_{1,n+1}$, $\ldots$, $P_{I_{n+1},n+1}$
of increased degrees.

A heuristic explanation of this construction is that the part
\[
\sum_{l\le T}\Bigl(Q_{1,\epsilon_1,\K,\alpha}-\epsilon_1Q_0\Bigr)\left(\frac{\log l}{\log T}+\delta_1(s)\right)l^{-s}
\]
of the function to be mollified, translatable by a positive $\ds$ to
\[
\sum_{l\le T}Q_{\K,\alpha}\left(\frac{\log l}{\log T}+\delta(s)\right)l^{-(s+\ds)},
\]
is expected to be better mollifiable than the non-translatable part
\[
\sum_{l\le T}\epsilon_1Q_0\left(\frac{\log l}{\log T}+\delta_1(s)\right)l^{-s}.
\]
But looking at the translatable part in the translated form we see
that the relative weight of the derivatives of $\zeta(s)$ is larger than the relative weight of $\zeta(s)$,
so we need more Feng's polynomials $P_{1}$, $\ldots$, $P_{I}$ of larger degrees for good mollification.

To obtain the value of $\kappa$ in Feng's method, consider
\[
F_{\epsilon_1,\K,\ds,R}(s)=G_{\epsilon_1,\K,\ds}(s)M_{\epsilon_1,\K,\ds,R}(s)
\]
with the mollifier
\[
\begin{split}
M_{\epsilon_1,\K,\ds,R}(s)&=\sum_{m\le T^{\theta}}\mu(m)P_{1,\epsilon_1,\K,\ds}\left(\frac{\log y/m}{\log y}\right)m^{-\left(s+\frac R{\log T}\right)}\\
&+\sum_{m\le T^{\theta_1}}\frac{\mu(m)}{m^{\frac{R}{\log
T}+s}}\Big(\sum_{p_1p_2|m}\frac{\log p_1\log p_2}{\log^2y_1}P_2\left(\frac{\log y_1/m}{\log y_1}\right)\\
&+\sum_{p_1p_2p_3|m}\frac{\log p_1\log p_2\log p_3}{\log^3y_1}P_3\left(\frac{\log y_1/m}{\log y_1}\right)+\cdots\\
&+\sum_{p_1p_2\cdots p_I|m}\frac{\log p_1\log p_2\cdots\log p_I}{\log^Iy_1}P_I\left(\frac{\log y_1/m}{\log y_1}\right)\Big),
\end{split}
\]
where $y=T^{\theta}$, $y_1=T^{\theta_1}$, $I\geq 2$ is an integer, $P_{1,\epsilon_1,\K,\ds}(x)$ is a real polynomial with
\[
P_{1,\epsilon_1,\K,\ds}(0)=0\qquad\text{and}\qquad P_{1,\epsilon_1,\K,\ds}(1)=1,
\]
$P_l (l=2,\ldots I)$ are real polynomials with
$P_l(0)=0$, $p_1,p_2,\ldots,p_I$ run over the prime numbers.

We then have $N_G(\RR)\le N_F(\RR)$ and applying the Littlewood lemma
for
\[
a=\frac12-\frac R{\log T},\quad|\ds|=\frac{|\alpha|}{\log T}=\frac{2\pi-\epsilon}{\log T}
\]
we arrive at the following principal inequality of the Levinson--Conrey method:
\begin{theorem}
\label{PreobLevCon}
Suppose that $\epsilon$, $R$ are fixed, $0<\epsilon<2\pi$, $R>0$,
$\K$ is a fixed large odd integer, $T$ goes to infinity,
\[
|\ds|=\frac{2\pi-\epsilon}{\log T},\quad\Re s=a=\frac12-\frac R{\log T}.
\]
Let $N_{00}(T,2T)$ be the number of zeros $s=\rho=\frac12+i\gamma$ of $\zeta(s)$
counted without multiplicity which are not zeros of $G_{\epsilon_1,\K,\ds}(s)$.
Then
\[
N_{00}(T,2T)\ge N(T,2T)\left(1-\frac2{R}\log I(R)+O\left(\frac1{\log T}\right)\right),
\]
where
\[
I(R)=\frac1T\inte_T^{2T}|F_{\epsilon_1,\K,\ds,R}(a+it)|\,dt.
\]
\end{theorem}

By this theorem, choosing
\[
\begin{split}
\theta&=\frac47-\varepsilon,\ \ \theta_1=\frac12-\varepsilon,\ \ R=1{.}3025,\ \ I=5,\\
P_1(x)&=x+0{.}138173\,x(1-x)-0{.}445606\,x(1-x)^2-4{.}039834\,x(1-x)^3\\
&+7{.}506942\,x(1-x)^4-3{.}239261\,x(1-x)^5,\\
P_2(x)&=-0{.}101269\,x+3{.}571698\,x^2-1{.}807283\,x^3-0{.}929884\,x^4,\\
P_3(x)&=1{.}334025\,x-3{.}018815\,x^2+1{.}133072\,x^3,\\
P_4(x)&=-0{.}546630\,x+0{.}372783\,x^2,\\
P_5(x)&=-1{.}029768\,x,\\
\alpha&=0{.}1,\qquad{\mathcal K}=5,\\
Q(x)&=Q_{1,\epsilon_1,{\mathcal K},\alpha}(x)=\ \text{Taylor expansion of order $\K$ of}\\
&\frac12-\frac12\tanh\left(\frac{\alpha}2\left(x-\frac12\right)\right)
+\left(\frac12-\frac12\tanh\left(\frac{\alpha}4\right)\right)(2\tilde{Q}(x)-1),\\
\tilde{Q}(x)&=1-0{.}6684\,x-1{.}0798\left(\frac{x^2}{2}-\frac{x^3}{3}\right)-5{.}0447\left(\frac{x^3}{3}-\frac{x^4}{2}+\frac{x^5}{5}\right),
\end{split}
\]
using the asymptotic formula given in~\cite[Theorem 2]{PreobFeng12} and making $\varepsilon\to0$ we obtain
\[
\kappa\ge0{.}412948.
\]

\textbf{4. Further remarks.}
Though not essential for our purposes, Lemma~\ref{PreobBernoulliLemma} has the limitation $|\alpha|<2\pi$.
This could be avoided if one proves that for $\alpha$ arbitrarily large
the analytic function given for $\Re s>1$ by
\[
g_{\alpha,T}(s)=-\frac12\sum_{l=1}^{\infty}\tanh\left(\frac{\alpha}2\left(\frac{\log l}{\log T}-\frac12\right)\right)l^{-s}
\]
obeys two types of symmetries:
\begin{enumerate}
\item For $\ds=\frac{\alpha}{\log T}$
\be\label{PreobS1}
2e^{\alpha/2}\left(\frac{\zeta(s+\ds)}2-g_{\alpha,T}(s+\ds)\right)=2\left(\frac{\zeta(s)}2+g_{\alpha,T}(s)\right).
\ee
\item For $s=\sigma+it$, $T\le t\le2T$, the function
\be\label{PreobS2}
H(s)g_{\alpha,T}(s)+\text{small perturbation}
\ee
is purely imaginary for $\Re s=\frac12$.
\end{enumerate}

To prove~\eqref{PreobS1} we note that
\[
\tanh\left(\frac u2\right)=\frac2{e^{-u}+1}-1
\]
and substitute this with
\[
u=\alpha\left(\frac{\log l}{\log T}-\frac12\right)
\]
into the left-hand side of~\eqref{PreobS1}, obtaining the Dirichlet series
\[
\sum_{l=1}^{\infty}\frac{2e^{\alpha/2}}
{\left(e^{-\alpha\left(\frac{\log l}{\log T}-\frac12\right)}+1\right)l^{s}e^{\alpha\frac{\log l}{\log T}}}.
\]
In the right-hand side of~\eqref{PreobS1} we use
\[
\tanh\left(\frac u2\right)=-\frac2{e^{u}+1}+1
\]
obtaining
\[
\sum_{l=1}^{\infty}\frac2
{\left(e^{\alpha\left(\frac{\log l}{\log T}-\frac12\right)}+1\right)l^{s}}.
\]
The two Dirichlet series are the same.

To prove~\eqref{PreobS2}, we note that
\[
\tanh\left(\frac u2\right)
\]
is an odd function of $u$, so for $u\in[-A,A]$
it can be uniformly approximated by finite sums of the odd powers of $u$.
Thus
\[
H(s)g_{\alpha,T}(s)
\]
is approximated by odd derivatives of the $\xi$ function with real coefficients,
which are purely imaginary on the critical line.

Details of the above argument are given in the Appendix.

\textbf{5. Appendix.}
\begin{lemma}[Analytic continuation of $g_{\alpha,T}(s)$]
\label{PreobAnContgLem}
For $s=\sigma+it$ with $\sigma>0$ and $0<t_0\le|t|\le2T$, where $t_0$ is fixed and $T\ge1$,
and for integer $N\ge T$ we have
\[
\begin{split}
g_{\alpha,T}(s)=-\frac12&\left(\sum_{n=1}^N\tanh\left(\frac{\alpha}2\left(\frac{\log n}{\log T}-\frac12\right)\right)n^{-s}\right.\\
&\left.{}-\frac{2T^{1-s}\log T}{\alpha(e^{\alpha/2}+1)(1-(1-s)(\log T)/\alpha)}F(1,1;2-(1-s)(\log T)/\alpha;(e^{\alpha/2}+1)^{-1})\right.\\
&\left.{}+\frac{T^{1-s}}{s-1}-\inte_T^{N+1/2}\tanh\left(\frac{\alpha}2\left(\frac{\log u}{\log T}-\frac12\right)\right)u^{-s}du\right.\\
&\left.{}+\frac{\alpha}{2\log T}\inte_{N+1/2}^{+\infty}\psi(u)\cosh^{-2}\left(\frac{\alpha}2\left(\frac{\log u}{\log T}-\frac12\right)\right)u^{-s-1}du\right.\\
&\left.{}-s\inte_{N+1/2}^{+\infty}\psi(u)\tanh\left(\frac{\alpha}2\left(\frac{\log u}{\log T}-\frac12\right)\right)u^{-s-1}du\right),
\end{split}
\]
where $F(a,b;c;z)$ is the hypergeometric function, and $\psi(x)=x-[x]-\frac12$.
\end{lemma}

\dokvo By the exact summation formula we have
\[
\begin{split}
&\sum_{N+1/2<n\le M+1/2}\tanh\left(\frac{\alpha}2\left(\frac{\log n}{\log T}-\frac12\right)\right)n^{-s}
=\inte_{N+1/2}^{M+1/2}\tanh\left(\frac{\alpha}2\left(\frac{\log u}{\log T}-\frac12\right)\right)u^{-s}du\\
{}+&\frac{\alpha}{2\log T}\inte_{N+1/2}^{M+1/2}\psi(u)\cosh^{-2}\left(\frac{\alpha}2\left(\frac{\log u}{\log T}-\frac12\right)\right)u^{-s-1}du\\
{}-&s\inte_{N+1/2}^{M+1/2}\psi(u)\tanh\left(\frac{\alpha}2\left(\frac{\log u}{\log T}-\frac12\right)\right)u^{-s-1}du.
\end{split}
\]
The first integral is convergent for $\sigma>1$ as $M\to+\infty$,
whereas the latter two integrals with the $\psi$ function converge absolutely for $\sigma>0$.
Denote them by $\Psi_1$ and $\Psi_2$.
Now for $\sigma>1$ we have the formula
\[
\begin{split}
g_{\alpha,T}(s)=-\frac12&\left(\sum_{n=1}^N\tanh\left(\frac{\alpha}2\left(\frac{\log n}{\log T}-\frac12\right)\right)n^{-s}\right.\\
&\left.{}+\inte_T^{+\infty}\tanh\left(\frac{\alpha}2\left(\frac{\log u}{\log T}-\frac12\right)\right)u^{-s}du
-\inte_T^{N+1/2}\tanh\left(\frac{\alpha}2\left(\frac{\log u}{\log T}-\frac12\right)\right)u^{-s}du\right.\\
&\left.{}+\frac{\alpha}{2\log T}\Psi_1-s\Psi_2\right),
\end{split}
\]
in which we consider
\[
\inte_T^{+\infty}\tanh\left(\frac{\alpha}2\left(\frac{\log u}{\log T}-\frac12\right)\right)u^{-s}du.
\]
We write the integrand as
\[
\tanh\left(\frac{\alpha}2\left(\frac{\log u}{\log T}-\frac12\right)\right)u^{-s}
=-\frac{2u^{-s}}{e^{-\alpha/2}u^{\alpha/\log T}+1}+u^{-s}.
\]
Integrating the latter term we get $\frac{T^{1-s}}{s-1}$,
while the former term gives
\[
\begin{split}
\inte_T^{+\infty}\frac{-2T^{-\frac s2+\frac12}}{\left(\frac u{\sqrt T}\right)^{\frac{\alpha}{\log T}}+1}
\left(\frac u{\sqrt T}\right)^{-s}d\left(\frac u{\sqrt T}\right)
&=-2T^{\frac{1-s}2}\inte_{\sqrt T}^{+\infty}\frac{x^{-s}dx}{1+x^{\alpha/\log T}}\\
&=\frac{-2T^{\frac{1-s}2}\log T}{\alpha}\inte_{e^{\alpha/2}}^{+\infty}\frac{v^{(1-s)(\log T)/\alpha-1}}{1+v}\,dv.
\end{split}
\]
Making the change of variables
\[
\begin{split}
w&=\frac1{v+1},\\
v&=\frac1w-1,\\
dv&=-\frac1{w^2}\,dw
\end{split}
\]
we get the integral
\[
\frac{-2T^{\frac{1-s}2}\log T}{\alpha}\inte_0^{(e^{\alpha/2}+1)^{-1}}w^{1-(1-s)(\log T)/\alpha-1}(1-w)^{(1-s)(\log T)/\alpha-1}dw
\]
that can be written as the incomplete beta function
\[
\frac{-2T^{\frac{1-s}2}\log T}{\alpha}B_{(e^{\alpha/2}+1)^{-1}}(1-(1-s)(\log T)/\alpha,(1-s)(\log T)/\alpha)
\]
which in turn can be expressed in terms of the hypergeometric function
\[
\begin{split}
&\frac{-2T^{\frac{1-s}2}\log T}{\alpha}(e^{\alpha/2}+1)^{(1-s)(\log T)/\alpha}\\
&\times F\left(1-(1-s)(\log T)/\alpha,1-(1-s)(\log T)/\alpha;2-(1-s)(\log T)/\alpha;(e^{\alpha/2}+1)^{-1}\right).
\end{split}
\]
Using the known linear transformation formula
\[
F(a,b;c;z)=(1-z)^{c-a-b}F(c-a,c-b;c;z)
\]
we get the term
\[
-\frac{2T^{1-s}\log T}{\alpha(e^{\alpha/2}+1)(1-(1-s)(\log T)/\alpha)}
F(1,1;2-(1-s)(\log T)/\alpha;(e^{\alpha/2}+1)^{-1})
\]
of the analytic continuation formula,
where the function
\[
F(1,1;2-(1-s)(\log T)/\alpha;(e^{\alpha/2}+1)^{-1})
\]
is analytic and bounded in $s$ for $|t|\ge t_0>0$ by the series representation.

\begin{lemma}[Approximate equation for $g_{\alpha,T}(s)$]
\label{PreobApproxEqgLem}
For $s=\sigma+it$ with $\sigma\ge\sigma_0>0$ and $0<t_0\le|t|\le2T$, where $\sigma_0$, $t_0$ are fixed and $T\ge1$,
we have
\[
\begin{split}
g_{\alpha,T}(s)=-\frac12&\left(\sum_{n=1}^T\tanh\left(\frac{\alpha}2\left(\frac{\log n}{\log T}-\frac12\right)\right)n^{-s}\right.\\
&\left.{}-\frac{2T^{1-s}\log T}{\alpha(e^{\alpha/2}+1)(1-(1-s)(\log T)/\alpha)}F(1,1;2-(1-s)(\log T)/\alpha;(e^{\alpha/2}+1)^{-1})\right.\\
&\left.{}+\frac{T^{1-s}}{s-1}+O\left(T^{-\sigma}\right)\right),
\end{split}
\]
where the constant in the $O$-term is absolute.
\end{lemma}

\dokvo In the analytic continuation formula of Lemma~\ref{PreobAnContgLem}
we use the standard uniform approximation
\[
\sum_{T<n\le N+1/2}\tanh\left(\frac{\alpha}2\left(\frac{\log n}{\log T}-\frac12\right)\right)n^{-s}
=\inte_T^{N+1/2}\tanh\left(\frac{\alpha}2\left(\frac{\log u}{\log T}-\frac12\right)\right)u^{-s}du+O\left(T^{-\sigma}\right)
\]
and make $N\to\infty$. $\qed$

We shall obtain approximations to $g_{\alpha,T}(s)$ by using the Fourier expansion
\[
\tanh\left(\frac{\alpha x}2\right)=\sum_{k=1}^Kb_k(\alpha)\sin(kx)+\hat{R}_{K,\alpha}(x)
\]
and the Taylor expansion
\[
\sin(kx)=\sum_{m=1}^M(-1)^{m-1}\frac{(kx)^{2m-1}}{(2m-1)!}+R_{k,M}(x),
\]
where
\[
\begin{split}
\hat{R}_{K,\alpha}(x)&=-2\inte_0^{\pi}\phi_{\alpha,x}(y)D_K(y)\,dy,\\
\phi_{\alpha,x}(y)&=\frac{\tanh\left(\frac{\alpha(x+y)}2\right)+\tanh\left(\frac{\alpha(x-y)}2\right)-2\tanh\left(\frac{\alpha x}2\right)}2,\\
R_{k,M}(x)&=\frac{(-1)^M(kx)^{2M+1}}{(2M)!}\inte_0^1(1-u)^{2M}\cos(kxu)\,du,
\end{split}
\]
and $D_K(y)$ is the Dirichlet kernel.
Explicitly, the coefficients $b_k(\alpha)$ are
\[
\begin{split}
b_k(\alpha)&=-\frac4{\pi}\inte_0^{\pi}\frac{e^{ikx}-e^{-ikx}}{(e^{\alpha x}+1)2i}\,dx\\
&=-\frac4{\pi\alpha 2i}\left(\inte_1^{e^{\alpha\pi}}\frac{v^{ik/\alpha-1}}{v+1}\,dv-\inte_1^{e^{\alpha\pi}}\frac{v^{-ik/\alpha-1}}{v+1}\,dv\right)\\
&=-\frac4{\pi\alpha}\Im\left(B_{1/2}\left(1-i\frac k{\alpha},i\frac k{\alpha}\right)-B_{(e^{\alpha\pi}+1)^{-1}}\left(1-i\frac k{\alpha},i\frac k{\alpha}\right)\right),
\end{split}
\]
where $B_x(a,b)$ is the incomplete beta function.

So we have
\be
\label{PreobOddPowApprox}
\begin{split}
&\sum_{l=1}^T\tanh\left(\frac{\alpha}2\left(\frac{\log l}{\log T}-\frac12\right)\right)l^{-s}\\
&=\sum_{k=1}^Kb_k(\alpha)\sum_{m=1}^M(-1)^{m-1}\frac{k^{2m-1}}{(2m-1)!}
\sum_{l\le T}\left(\frac{\log l}{\log T}-\frac12\right)^{2m-1}l^{-s}\\
&+\sum_{l\le T}{\mathcal R}_{K,M,\alpha}\left(\frac{\log l}{\log T}-\frac12\right)l^{-s},
\end{split}
\ee
where
\[
{\mathcal R}_{K,M,\alpha}(x)=\hat{R}_{K,\alpha}(x)+\sum_{k=1}^Kb_k(\alpha)R_{k,M}(x).
\]
We multiply the polynomial
\[
\sum_{k=1}^Kb_k(\alpha)\sum_{m=1}^M(-1)^{m-1}\frac{k^{2m-1}}{(2m-1)!}
\left(\frac{\log l}{\log T}-\frac12\right)^{2m-1}
\]
appearing in the right-hand side of~\eqref{PreobOddPowApprox} by $-\frac12$ and denote it
$q\left(\frac{\log l}{\log T}\right)$.

\begin{lemma}[Approximation to $H(s)g_{\alpha,T}(s)$ by a sum of the odd derivatives of the $\xi$ function]
\label{PreobApproxHgLem}
For $s=\sigma+it$ in the rectangle $\frac13\le\sigma\le A$, $T\le t\le2T$, with $A\ge3$ and $T\ge2A$,
we have
\[
2H(s)\left(\sum_{l\le T}q\left(\frac{\log l}{\log T}+\delta(s)\right)l^{-s}+O(T^{-\frac14})\right)
=\sum_{m=1}^Mg_{2m-1}\xi^{(2m-1)}(s),
\]
where
\[
\delta(s)=\frac{\log(2\pi T/s)}{2\log T}\ll\frac1{\log T},
\]
and $g_{2m-1}$ are real numbers.
\end{lemma}

\dokvo See~\cite[Chapter 18]{PreobIw14}.

Now in Section~3 we can substitute
\[
Q_{1,\epsilon_1,{\mathcal K},\alpha}(x)=\frac12+\epsilon_1Q_0(x)+q(x)
\]
(with $2M-1$ in place of ${\mathcal K}$),
and
\[
Q_{{\mathcal K},\alpha}(x)=\frac12-q(x)
\]
with $K$, $M$, $|\alpha|$ arbitrarily large.

\bigskip

\end{document}